\documentclass[a4paper,11pt]{article}
\usepackage{amsfonts}
\usepackage{amsmath}
\usepackage{amssymb}
\usepackage{amstext}
\usepackage{amsthm}

\usepackage{xspace}
\usepackage[dvips]{graphicx}

\newcommand{\C}{\mathbb C}
\newcommand{\R}{\mathbb R}

\newcommand{\E}{{\mathcal E}}
\newcommand{\Ch}{{\mathcal C}}
\newcommand{\eps}{\varepsilon}

\newcommand{\set}[1]{\left\{#1\right\}}

\newtheorem{theorem}{Theorem}[section]
\newtheorem{lemma}[theorem]{Lemma}
\newtheorem{prop}[theorem]{Proposition}
\newtheorem{cor}[theorem]{Corollary}

\theoremstyle{definition}
\newtheorem{definition}[theorem]{Definition}
\newtheorem{rem}[theorem]{Remark}
\theoremstyle{remark}

\numberwithin{equation}{section}

\begin{document}

\title{Solitons for the nonlinear Klein-Gordon equation}
\author{J. Bellazzini \and V. Benci \and C. Bonanno \and A.M. Micheletti \\[0.3cm]
Dipartimento di Matematica Applicata "U. Dini"\\
Universit\`a di Pisa, Via F. Buonarroti 1/c, Pisa, ITALY. \\
email:\ j.bellazzini@ing.unipi.it,\ benci@dma.unipi.it,\\
bonanno@mail.dm.unipi.it,\,
a.micheletti@dma.unipi.it.}
\date{}
\maketitle

\begin{abstract}
In this paper we study existence and orbital stability for
solitary waves of the nonlinear Klein-Gordon equation. The energy
of these solutions travels as a localized packet, hence they are a
particular type of solitons. In particular we are interested in
sufficient conditions on the potential for the existence of
solitons. Our proof is based on the study of the ratio
energy/charge of a function, which turns out to be a useful
approach for many field equations.
\end{abstract}

\section{Introduction} \label{intro}

We are interested in studying the existence of a particular class
of solutions for the nonlinear wave equation: the soliton
solutions. By \textit{solitary wave} we mean a solution of a field
equation whose energy travels as a localized packet; if a solitary
wave exhibits orbital stability is called \textit{soliton}. In
this respect solitons have a particle-like behavior and they occur
in many questions of mathematical physics, such as classical and
quantum field theory, nonlinear optics, fluid mechanics, plasma
physics (see e.g. \cite{Witham},\cite{Rajaraman}).

Today, we know (at least) three mechanisms which might produce solitary
waves and solitons:
\begin{itemize}
\item  (i) Complete integrability; e.g. Korteweg-de Vries equation
\begin{equation}
u_{t}+u_{xxx}+6uu_{x}=0  \tag{KdV}  \label{KdV}
\end{equation}

\item  (ii) Topological constrains: e.g. Sine-Gordon equation
\begin{equation}
u_{tt}-u_{xx}+\sin u=0  \tag{SNG}  \label{SNG}
\end{equation}

\item  (iii) Ratio energy/charge: e.g. the nonlinear Klein-Gordon equation
\begin{equation}
\psi _{tt}-\Delta \psi +W^{\prime }(\psi )=0  \tag{NKG}
\label{eq}
\end{equation}
where $\psi (t,x)\in H^{1}(\R \times \R^{N},\C)$ with $N\geq 2$,
and $W:\C \rightarrow \R $ with $W(\psi )=F(|\psi |)$ for some
$C^{2}$ function $F:\R \rightarrow \R$. This kind of solitons are
called \emph{hylomorphic solitons} in \cite{hylo}.
\end{itemize}

The study of solitons for equation (\ref{eq}) has a very long
history starting with the pioneering paper of Rosen
\cite{rosen68}. Coleman \cite {coleman78} and Strauss
\cite{strauss} gave the first rigorous proofs of existence of
solutions of the type (\ref{standing-wave-2}) for some particular
$W^{\prime}s$, and later necessary and sufficient existence
conditions have been found by Berestycki and Lions
\cite{Beres-Lions}.

The first orbital stability results for (\ref{eq}) are due to
Shatah; in \cite{shatah} a necessary and sufficient condition for
orbital stability is given. We will discuss this result in Section
\ref{orbit-st}. However this condition is difficult to be verified
in concrete situations due to the fact that one should know
explicit properties of the solitary wave. See also \cite{gss87}
for a generalization of methods in \cite{shatah}.

The aim of this paper is to establish sufficient conditions on the
function $W$ which guarantee the existence of stable solitary
waves. Moreover these conditions are in some sense also necessary
(see Remark \ref{condizioni-necessarie} below).

The study of solitary waves of equation (\ref{eq}) can be reduced
to that of standing waves. A \textit{standing wave} is a finite
energy solution of (\ref{eq}) having the following form
\begin{equation}  \label{standing-wave-2}
\psi(t,x)=u(x)e^{-i\omega t}, \ u \geq 0,\ \omega \in \R .
\end{equation}
Substituting (\ref{standing-wave-2}) in (\ref{eq}), we get
\begin{equation}
-\Delta u+W^{\prime }(u)=\omega^{2}u  \label{static}
\end{equation}
Since equation (\ref{eq}) is obtained by studying critical points
of a Lagrangian with density $L$ (see (\ref{lag})) which is
invariant for the Poincar\'e group, we can obtain solutions $\psi
'(t,x)$ just making a Lorentz transformation of a given solution
$\psi(t,x)$. Namely, if we take the velocity $\mathbf{v}=(v,0,0),$
$ \left| v\right| <1$, and set
\begin{equation*}
t^{\prime }=\gamma \left( t-vx_{1}\right) ,\ x_{1}^{\prime
}=\gamma \left( x_{1}-vt\right) ,\ x_{2}^{\prime
}=x_{2},\ x_{3}^{\prime }=x_{3}\ \ \text{with}\ \ \gamma
=\frac{1}{\sqrt{1-v^{2}}}
\end{equation*}
it turns out that
\begin{equation*}
\psi '(t,x)=\psi (t^{\prime },x^{\prime })
\end{equation*}
is a solution of (\ref{eq}).

In particular given a standing wave $\psi (t,x)=u(x)e^{-i\omega
t},$ the function $\psi _{\mathbf{v}}(t,x):=\psi
(t^{\prime},x^{\prime })$ is a solitary wave which travels with
velocity $\mathbf{v}$. Thus, if $u(x)=u(x_{1},x_{2},x_{3})$ is any
solution of equation (\ref{static}), then
\begin{equation*}
\psi _{\mathbf{v}}(t,x_{1},x_{2},x_{3})=u\left( \gamma \left(
x_{1}-vt\right) ,x_{2},x_{3}\right) e^{i(\mathbf{k\cdot x}-\bar \omega t)},
\end{equation*}
is a solution of (\ref{eq}) provided that
\begin{equation*}
\bar \omega =\gamma \omega\quad  \text{and} \quad \mathbf{k}=\gamma \omega
\mathbf{v}
\end{equation*}
Hence to study the existence of solitons we can restrict to
consider standing waves.

We consider the following hypotheses on the potential $W$:
\begin{equation}
W(0)=W^{\prime }(0)=0,\,W^{\prime \prime }(0)=1  \tag{$H_{0}$}  \label{H0}
\end{equation}
\begin{equation}
\alpha _{0}:=\inf\limits_{z\in \C }\ \frac{W(z)}{\frac{1}{2}|z|^{2}}<1
\tag{$H_{1}$}  \label{H1}
\end{equation}
\begin{equation}
W(z)\geq 0  \tag{$H_{2}$}  \label{H2}
\end{equation}

Our main result is the following

\begin{theorem} \label{main-theorem}
If (\ref{H0}),(\ref{H1}) and (\ref{H2}) hold then equation
(\ref{eq}) admits solitons.
\end{theorem}

\begin{rem}
We now make some comments on the assumptions of the theorem above.
The assumption $W^{\prime \prime }(0)=1$ is a normalization
condition; making a change of the space variable, (\ref{H0}) is
satisfied provided that $W^{\prime \prime }(0)>0$. We remark that
if $W^{\prime \prime}(0)<0$ there are no standing waves, and if
$W^{\prime \prime}(0)=0$ standing waves exist under restrictive
assumptions on the growth of $W$ (see \cite{Beres-Lions}).
Condition (\ref{H1}) is the crucial assumption and it is typical
for solitons of type (iii). Assumption (\ref{H2}) is natural for
many physical models described by (\ref{eq}); however it is not
necessary to have orbital stability (see also Remark
\ref{condizioni-necessarie} below).
\end{rem}

\section{Statement of the main results} \label{main-results}

\subsection{The wave equations as dynamical systems} \label{wave-ds}

Let us consider the Lagrangian density
\begin{equation}  \label{lag}
L(\psi)=\frac{1}{2}\left( \left| \frac{\partial \psi }{\partial t}\right|
^{2}-|\nabla \psi |^{2}\right) -W(\psi ).
\end{equation}
Its critical points satisfy equation (\ref{eq}), with
\begin{equation*}
W^{\prime }(\psi )=F^{\prime }(\left| \psi \right| )\frac{\psi }{\left| \psi \right| }.
\end{equation*}
Moreover $L(\psi)$ is invariant for the action of the Poincar\'{e}
group and for the gauge action of the group $S^{1}$ given by
\begin{equation} \label{gauge}
\psi \mapsto e^{i\theta }\psi \qquad \theta \in \R
\end{equation}
and can be considered as the simplest nonlinear Lagrangian density with these invariance properties. In the following sections we will see that equation (\ref{eq}), with assumptions (\ref{H0})-(\ref{H2}) on $W$, produces a very rich model in which there are solitary waves which behave as relativistic particles.

We remark that if $W^{\prime }(\psi)$ is linear, namely $W^{\prime }(\psi
)=\Omega ^{2}\psi$, then equation (\ref{eq}) reduces to the Klein-Gordon
equation. Among the solutions of the Klein-Gordon equations there are the
\textit{wave packets} which behave as solitary waves but disperse in space
as time goes on. On the contrary, if $W$ has a nonlinear suitable
component the wave packets do not disperse, hence one can find solitons.

We now reconsider the problem from the point of view of dynamical systems.
Let us consider the space $X= H^1(\R^N,\C)\times L^2(\R^N,\C)$. To a function $\psi(t,x) \in H^1(\R \times \R^N,\C )$ we can associate a vector $\mathbf{\Psi}
=(\psi,\psi_t)$ which is in $X$ for any $t\in \R$. Equation (\ref{eq}) can be
written in the form
\begin{equation}  \label{eq-hamil}
\frac{\partial \mathbf{\Psi}}{\partial t} = A \mathbf{\Psi} - \mathbf{
W^{\prime}}(\mathbf{\Psi})
\end{equation}
for the matrix
\begin{equation*}
A = \left(
\begin{array}{cc}
0 & 1 \\[0.2cm]
\triangle & 0
\end{array}
\right)
\end{equation*}
and vector function
\begin{equation*}
\mathbf{W^{\prime}}(\mathbf{\Psi}) = \left(
\begin{array}{c}
0 \\[0.2cm]
W^{\prime}(\psi)
\end{array}
\right).
\end{equation*}
Together with an initial condition $\mathbf{\Psi}(0,x)$, equation
(\ref {eq-hamil}) is a Cauchy problem, and its solutions define a
flow $U$ on the phase space $X$, that is a map $U:\R \times X \to
X$ given by
\begin{equation}  \label{flusso}
\mathbf{\Psi}(t,x)=U(t,\mathbf{\Psi}(0,x))
\end{equation}
Global existence for the flow $U$ is guaranteed by assumption (\ref{H2}).

Noether's theorem states that any invariance for a one-parameter group of
the Lagrangian implies the existence of an integral of motion, namely of a
quantity on solutions which is preserved with time (see e.g. \cite{Gelfand}). Thus equation (\ref{eq-hamil}) has ten integrals: energy,
momentum, angular momentum and ergocenter velocity. Moreover, another
integral is given by the gauge invariance (\ref{gauge}): charge. We now
briefly introduce these integrals, for details we refer to \cite{Gelfand}.
In this paper we are only interested in energy and charge.

Energy, by definition, is the quantity which is preserved by the time
invariance of the Lagrangian. In our case it has the following form
\begin{equation}  \label{energy}
\E(\mathbf{\Psi})=\int \left[ \frac{1}{2}\left| \psi_t \right|^{2} + \frac{1
}{2} \left| \nabla \psi \right|^{2} + W(\psi) \right]\, dx
\end{equation}

Momentum and angular momentum are the quantity preserved respectively by the space translations and space rotations invariance of the Lagrangian. Their expressions are given by
\begin{equation*}
\vec{P}(\mathbf{\Psi}) = -\mathrm{Re} \int \psi_t \overline{\nabla \psi }\;
dx
\end{equation*}
\begin{equation*}
\vec{M}(\mathbf{\Psi}) = \mathrm{Re} \int \left( \mathbf{x}\times \nabla
\psi \right)\, \overline{\psi_t}\ dx
\end{equation*}

If we take the Lagrangian (\ref{lag}), the quantity preserved by the Lorentz
transformation, by standard computations, is the derivative of the
ergocenter $\vec{Q}$ which satisfies
\begin{equation*}
\dot{\vec{Q}}(\mathbf{\Psi})= \frac{\vec{P}(\mathbf{\Psi})}{\E(\mathbf{\Psi})}
\end{equation*}
Then, the three components of $\dot{\vec{Q}}$ are other three integrals of
motion.

Finally, we consider the charge, which is the quantity preserved by the
trivial gauge action (\ref{gauge}). The charge for a general Lagrangian $L$
has the following expression (see e.g. \cite{sammomme})
\begin{equation*}
\Ch(\mathbf{\Psi}) = \mathrm{Im} \int \frac{\partial L}{\partial \psi_t} \,
\overline{\psi }\ dx
\end{equation*}
In our case we get
\begin{equation*}
\Ch(\mathbf{\Psi}) =\mathrm{Im} \int \psi_t\, \overline{\psi }\ dx
\end{equation*}
Without loss of generality we consider solutions with negative charge. The
same results can be obtained in the same way also in the case of positive
charge.

As stated in the introduction, a particular class of solutions of (\ref{eq}) (and of (\ref{eq-hamil})) is given by standing waves, namely solutions of the form (\ref{standing-wave-2}) for which the couple $(u,\omega)$ satisfies equation (\ref{static}).

Moreover if $u\in H^1(\R^N,\R^{+})$ is a solution of (\ref{static}) for any given $\omega\in \R$, we have a standing wave for (\ref{eq}). For a standing wave $\psi_{0}$ the relative element in the space $X$ is given by $\mathbf{\Psi_0}=(u(x)e^{-i\omega t},-i\omega u(x)e^{-i\omega t})$, and energy and charge
are given by
\begin{equation}  \label{enerstand}
\E(\mathbf{\Psi_0}) = \int \left( \frac{1}{2} |\nabla u|^2+ \frac 1 2
\omega^2 u^2 +W(u) \right)dx,
\end{equation}
\begin{equation}  \label{carstand}
\Ch(\mathbf{\Psi_0})= - \int\, \omega\, u^2\ dx.
\end{equation}
Since we consider only the case of negative charge, it holds $\omega \in \R^+$.

We can identify the standing waves as a subset of
\begin{equation}  \label{set-u-omega}
Y:= \left\{(u,\omega) \in H^1(\R^N,\R^+)\times \R^+\right\}
\end{equation}
through the embedding of $Y$ into $X$ given by
\begin{equation}  \label{imm-y-x}
Y \ni (u,\omega) \mapsto (u(x)e^{-i\omega t}, -i \omega u(x)e^{-i\omega t})
\in X \quad \forall\, t \in \R
\end{equation}
and define two functionals on $Y$, energy
\begin{equation}  \label{energy-y}
E(u,\omega) :=\int \left( \frac{1}{2} |\nabla u|^2+ \frac 1 2 \omega^2 u^2
+W(u) \right)dx,
\end{equation}
and charge
\begin{equation}  \label{charge-y}
C(u,\omega)= - \int\, \omega\, u^2\ dx.
\end{equation}
Notice that energy $\E$ and charge $\Ch$ of a standing wave $\mathbf{\Psi_0}$ ((\ref{enerstand}) and (\ref{carstand})) coincide with energy $E$ and charge $C$ of the couple $(u,\omega) \in Y$ corresponding to $\mathbf{\Psi_0}$.

Using these definitions, standing waves $\mathbf{\Psi_0}$ will be found as points
of constrained minima of the energy function $E(u,\omega)$ on the manifold
\begin{equation}  \label{constrain}
M_C=\left\{ (u,\omega)\in Y: C(u,\omega)=C\right\}
\end{equation}
of $(u,\omega)$ with fixed charge $C\not= 0$.

In the following, an important quantity is the ratio between energy and absolute value of charge, hence we introduce
the notation $\Lambda(u,\omega)$ for the functional on $Y$ defined by
\begin{equation}  \label{rapporto}
\Lambda(u,\omega) := \frac{E(u,\omega)}{|C(u,\omega)|} = \frac 1 2 \, \omega
+ \frac{1}{2\omega} \, \alpha(u)
\end{equation}
where
\begin{equation}  \label{alpha}
\alpha(u):= \frac{\int \left( \frac{1}{2} |\nabla u|^2 + W(u) \right)dx}{
\int \frac 1 2 \, u^2 \ dx}
\end{equation}
The study of the functional $\Lambda$ turns out to be useful for a general approach to many field equations. For a discussion of these subjects we refer to \cite{hylo}.

\begin{prop}
\label{alpha-0} Let $\alpha_0$ be defined as in (\ref{H1}). Then
\begin{equation*}
\alpha_0 = \inf\limits_{u\in H^1} \ \alpha(u)
\end{equation*}
\end{prop}

\begin{proof}
By (\ref{H1}), for all $u\in H^{1}$ it holds
$$
\alpha(u) \ge \frac{\int  W(u)\, dx}{\int \frac 1 2 \, u^2 \ dx}\, \ge \frac{\alpha_{0} \int \frac 1 2 \, u^2 \ dx}{\int \frac 1 2 \, u^2 \ dx} = \alpha_{0}
$$
hence one inequality is proved. For the opposite inequality, by (\ref{H1}) for any $\eps >0$ there exists $z_{0}\in \C$ such that
$$
\frac{W(z_{0})}{\frac 1 2 |z_{0}|^{2}} < \alpha_{0} + \frac \eps 2
$$
We remark that $W$ depends only on the modulus of $z$, and denote $r_{0}=|z_{0}|$. Let us choose a family of functions $u_{R}(x) \in H^{1}$, with $R\ge 0$, by
\begin{equation}\label{frittate}
u_{R}(x)= \left\{
\begin{array}{cl}
r_{0} & \mbox{ if $|x| \le R$} \\[0.3cm]
0 & \mbox{ if $|x| \ge R+1$} \\[0.3cm]
r_{0} (1+R-|x|) & \mbox{ if $R\le |x| \le R+1$}
\end{array} \right.
\end{equation}
Moreover, let $\mu$ be the $N$-dimensional Lebesgue measure, $B_{R}=\set{|x|\le R}$ and $B_{R}^{R+1}=\set{R\le |x|\le R+1}$, then we need the basic result
\begin{equation} \label{misura-palle}
\lim\limits_{R \to \infty}\ \frac{\mu(B_{R}^{R+1})}{\mu(B_{R})} \, =0
\end{equation}
Hence, by (\ref{frittate}) and (\ref{misura-palle}) for $R$ big enough it holds
$$
\alpha(u_{R}) \le \frac{\int_{B_{R}^{R+1}} \left( \frac{1}{2} |\nabla u_{R}|^2 + W(u_{R}) \right) dx + \int_{B_{R}} W(u_{R}) dx}{\int_{B_{R}} \frac 1 2 \, u_{R}^2 \ dx } \le
$$
$$
\le \frac{\mu(B_{R}^{R+1}) (\frac 1 2 r_{0}^{2} + \max_{0\le |z|\le r_{0}} W(z))}{\frac 1 2 r_{0}^{2}\, \mu(B_{R})} + \frac{W(z_{0})}{\frac 1 2 |z_{0}|^{2}} < \alpha_{0} + \eps
$$
Hence the proposition in proved.
\end{proof}

\subsection{Orbital stability} \label{orbit-st}

An important property of the orbits of a dynamical system is their stability
with respect to small perturbations. In the following we consider on $X$ the
distance induced by the usual product topology of $H^{1}(\R^{N},\C) \times
L^{2}(\R^{N},\C)$.

\begin{definition} \label{stabilita-insiemi}
Let $\Gamma$ be a subset of $X$ which is invariant
under the evolution flow $U:\R \times X \to X$ of (\ref{eq-hamil}) (see (\ref
{flusso})). The set $\Gamma$ is called \emph{stable} if for any $\varepsilon
>0$ there exists $\delta >0$ such that $d(\Gamma, \mathbf{\Psi}(0,x))<
\delta $ implies that $d(\Gamma, U(t,\mathbf{\Psi}(0,x))< \varepsilon$ for
all $t\in \R$.
\end{definition}

A standing wave is then defined to be orbitally stable if its
orbit $U$ is a stable set under the flow in the sense of
Definition \ref{stabilita-insiemi}. Moreover, since equation
(\ref{eq}) is invariant under space and time translations, we
consider also all the possible translates of its orbit.

\begin{definition}
\label{stabilita-orbitale} A standing wave $\psi_0(t,x)=u(x)e^{-i\omega t}$
is called \emph{orbitally stable} if the set
{\small
\begin{equation}  \label{gamma-orbita}
\Gamma(u,\omega):= \left\{ (u(x+y)\, e^{-i(\omega t -\theta)}, -i\omega
u(x+y)\, e^{-i(\omega t -\theta)}) \ :\ y\in \R^{N},\, \theta \in \R \right\}
\subset X
\end{equation}}
is stable for the evolution flow $U$.
\end{definition}

In \cite{shatah} and \cite{gss87}, the orbital stability of a standing wave is studied with respect to its frequency $\omega$. Namely, the authors consider the real function
$$
\omega \mapsto d(\omega):= E(u_{\omega},\omega) + \omega C(u_{\omega},\omega)
$$
where $E$ and $C$ are the energy and the charge on $Y$, and
$u_{\omega}$ is the ground state solution of (\ref{static}) for a
fixed $\omega$. They prove that a necessary and sufficient
condition for the orbital stability of a standing wave
$\psi(t,x)=u(x) e^{-i\omega_{0} t}$ is the convexity of the map
$d(\omega)$ in $\omega_{0}$. This condition is difficult to be
verified since $d(\omega)$ cannot be computed explicitly. However
in some particular cases the properties of $d(\omega)$ can be
investigated. For example in \cite{ss85} the authors give exact
ranges of the frequency $\omega$ for which they obtain stability
and instability of the respective standing waves for the nonlinear
wave equation with $W(\psi) = \frac 1 2 |\psi|^{2} - \frac 1 p
|\psi|^{p}$. In particular stability occurs for $2< p < 2 + 4/N$.
Unfortunately, for a general $W$, it is very difficult to prove or
disprove the convexity of $d(\omega)$.

The aim of this paper is to overcome this difficulty and to establish the stability of standing waves looking directly at $W$ itself.

As stated above, standing waves will be found as constrained points of
minimum $(u,\omega)$ for the energy functional $E$ on the space $Y\subset X$
with fixed charge $C$. We prove stability for these standing waves by
proving stability with respect to the flow $U$ of the \emph{ground state set}
$\Gamma_{C}$, that is defined for a given charge $C$ as the set
\begin{equation}  \label{gamma-totale}
\Gamma_C:= \bigcup_{E(u, \omega)=\min_{_{M_C}} E} \ \Gamma(u,\omega)
\end{equation}

To prove stability of the ground state set we use the Lyapunov stability
method.

\begin{definition}
\label{lyapunov} A function $V:X \rightarrow \R$ is a Lyapunov function for
the set $\Gamma_C$ if

\begin{itemize}
\item  $V\geq 0$ and $V(\mathbf{\Psi})=0$ if and only if the function $\mathbf{\Psi}$ is of
the form $(u(x) e^{-i\omega t},-i\omega u(x) e^{-i\omega t})$ with $
(u,\omega) \in \Gamma_C$;

\item  $\frac{d}{dt}V(\mathbf{\Psi}(t,x)) \leq 0$ for any $t$, where $
\mathbf{\Psi}(t,x)=U(t,\mathbf{\Psi}(0,x))$ is solution of (\ref{eq-hamil});

\item  if $V(\mathbf{\Psi_{n}}) \rightarrow 0$ then $d(\mathbf{\Psi_{n}}
,\Gamma_{C}) \rightarrow 0$.
\end{itemize}
\end{definition}

\begin{prop}[Lyapunov stability method]
If there exists a a Lyapunov function for the set $\Gamma_C$, then $\Gamma_C$
is stable.
\end{prop}

\subsection{Results} \label{results}

We now state our results. The proofs are given in Section \ref{proof-main}.

Our main result is Theorem \ref{main-theorem}. It follows from some preliminary results which have some interest in themselves. We first prove the existence of a standing wave using a new proof based on the following variational principle. Recall (\ref{set-u-omega}) and (\ref{constrain}) for the definition of $Y$ and $M_{C}$.

\begin{theorem} \label{thm-critical}
A standing wave $\psi_0(t,x)=u(x)e^{-i\omega t}$ is a solution of (\ref{eq}) if and only if $(u,\omega) \in Y$ is a critical point of $E(u,\omega)$ constrained on the manifold $M_C$  with $C=C(u,\omega)$.
\end{theorem}

\begin{lemma}
\label{lemma-existence} If there exists $(\bar u, \bar \omega) \in Y$ such
that $\Lambda(\bar u, \bar \omega) <1$, there exist standing waves solutions
of (\ref{eq}) obtained as points of minima of the energy $E(u,\omega)$
constrained to the manifold $M_{\bar C}$ with $\bar C=C(\bar u,\bar \omega)$.
\end{lemma}

By (\ref{H1}) and Proposition \ref{alpha-0}, it follows that there exists $
\bar u$ such that $\alpha(\bar u)<1$. Moreover, by definition of $\Lambda$
(see (\ref{rapporto})), it follows that
\begin{equation*}
\inf\limits_{\omega \in \R^+} \, \Lambda(\bar u,\omega) = \sqrt{\alpha(\bar u)}
<1
\end{equation*}
Hence there exists $(\bar u, \bar \omega)$ such that $\Lambda(\bar u, \bar
\omega)<1$. Hence Lemma \ref{lemma-existence} implies the existence of
standing waves solutions of (\ref{eq}). Then we prove the orbital stability
of the standing waves.

\begin{theorem} \label{thm-stability}
If $(u,\omega)$ is an isolated point of local
minimum of the functional $E(u,\omega)$ constrained to the manifold $M_C$
with $C=C(u,\omega)$, then $\psi_0(t,x)=u(x)e^{-i\omega t}$ is an orbitally
stable standing wave.
\end{theorem}

Theorem \ref{main-theorem} is a pure existence result and gives no
information on the charge and frequency of the standing waves. Using Lemma
\ref{lemma-existence}, we obtain results for the charge. Properties of the
standing waves for different charges with respect to the form of the
potential $W$ and relative frequencies are discussed in details in \cite
{hylo}.

\begin{cor}
\label{123-carica} If (\ref{H0}),(\ref{H1}),(\ref{H2}) hold then there
exists $C_0>0$ such that for any $|C| \in (C_0,\infty)$ there exists an
orbitally stable standing wave for equation (\ref{eq}) with charge $C$.
\end{cor}

Existence of stable standing waves for all charges can be obtained by a
stronger version of assumption (\ref{H1}). Let us consider the condition

\begin{equation}  \label{H1-prime}
W(z) - \frac 1 2 |z|^2 \leq -|z|^{2+\varepsilon} \text{ with } 0<\varepsilon<
\frac{4}{N} \text{ for all $z\in \C$ with $|z|$ small.}  \tag{$H'_1$}
\end{equation}

\begin{cor}
\label{123-carica-tutte} If (\ref{H0}),(\ref{H1-prime}),(\ref{H2}) hold then
for any $C\not= 0$ there exists an orbitally stable standing wave for
equation (\ref{eq}) with charge $C$.
\end{cor}

\bigskip

\begin{rem} \label{condizioni-necessarie}
If condition (\ref{H3}) below holds, it is possible to prove that assumption (\ref{H2}) is not only sufficient but also necessary in order to get the existence of an absolute minimum of the energy $E(u,\omega)$ constrained to the manifold $M_{\bar{C}}.$ Thus, if
this assumption is violated, it is still possible to have orbitally stable
solitary waves as points of local minimum, but they can be destroyed by a perturbation which send them out of the basin of attraction.
\end{rem}

\section{Proofs} \label{proof-main}

In this section we consider only the case $C<0$ and $\omega \in \R^{+}$. For $C>0$ and $\omega \in \R^{-}$ the proofs are identical.

Moreover we give proofs of our results and in particular of Theorem \ref{main-theorem} which hold under an additional assumption. By (\ref{H0}),(\ref{H1}),(\ref{H2}), we can write
\begin{equation}  \label{R-u}
W(\psi)= \frac 1 2 |\psi|^2 + R(|\psi|)
\end{equation}
where $R: \R^+ \to \R$ is a $C^2$ function satisfying
\begin{equation*}
R(0)=R^{\prime}(0)=R^{\prime\prime}(0)=0, \quad R(s) \ge - \frac 1 2 s^{2} \ \ \forall \, s \in \R^{+}
\end{equation*}
and there exists $s_0\in \R^+$ such that $R(s_0)<0$. We assume that for all $s\in \R^{+}$
\begin{equation} \label{H3}
|R''(s)|\leq c_{1}s^{p-2}+c_{2}s^{q-2}\ \text{ with }c_{1},c_{2}>0 \text{
for }2<p\le q<2^{\ast }  \tag{$H_{3}$}
\end{equation}
This doesn't change the generality of our results. Indeed if (\ref{H3}) is violated for the potential $W$, we can define a potential $\tilde W$ which satisfies all assumptions (\ref{H0})-(\ref{H3}) and for which there exists $\bar s >0$ such that
\begin{equation} \label{uguali}
\tilde W(z) = W(z) \qquad \forall\, |z|< \bar s,
\end{equation}
and the functions $\tilde R(|z|) = \tilde W(z) - \frac 1 2 |z|^{2}$ satisfies
\begin{equation} \label{crescenza}
\tilde R'(s) \ge 0 \qquad \forall\, s \ge \bar s.
\end{equation}
Hence, by our results it follows that there exists a charge $C$ and a stable standing wave $\psi (t,x)= u(x) e^{-i  \omega t}$ for which
\begin{equation} \label{static-tilde}
-\Delta  u + \tilde W'(u) = \omega^{2} u
\end{equation}
and $(u,\omega)$ is a point of minimum for the energy $\tilde E$ on the manifold $M_{C}$, where $\tilde E$ is the energy of equation (\ref{energy-y}) with $\tilde W$ instead of $W$. We have the following result for $\psi$.

\begin{lemma} \label{principio-massimo}
Let $\tilde W$ satisfy (\ref{H0})-(\ref{H3}) and (\ref{crescenza}), then
$$
\| \psi(t,\cdot) \|_{L^{\infty}(\R^{N})} \le \bar s \qquad \forall\, t\in \R
$$
\end{lemma}

\begin{proof} Let $u$ be a solution of (\ref{static-tilde}) and set $u = \bar s + v$. It is sufficient to prove that $v\le 0$. Let $\Omega := \set{x\ :\ v(x) \ge 0}$. By (\ref{static-tilde}) we have that
$$
\begin{array}{ll}
-\Delta  v + \tilde W'(\bar s +v) = \omega^{2} (\bar s +v) & \mbox{in $\Omega$} \\[0.3cm]
v = 0 & \mbox{on $\partial \Omega$}
\end{array}
$$
Multiplying both sides of the above equation by $v$ and integrating in $\Omega$, we
get
$$
\begin{array}{cl}
0 & =  \int_{\Omega} \, \left[ |\nabla v|^{2} + \tilde W'(\bar s +v)v - \omega^{2} (\bar s + v)v \right] \, dx \\[0.3cm]
& =  \int_{\Omega} \, \left[ |\nabla v|^{2} + \tilde R'(\bar s +v)v + (1 - \omega^{2}) (\bar s + v)v \right] \, dx \\[0.3cm]
& \ge  \int_{\Omega} \, \left[ |\nabla v|^{2} + \tilde R'(\bar s +v)v  \right] \, dx \\[0.3cm]
& \ge  \int_{\Omega} \, |\nabla v|^{2} \, dx
\end{array}
$$
where we have used (\ref{crescenza}) and $\omega <1$, which follows from Remark \ref{remark-conv}. From this it follows that $v=0$ in $\Omega$.
\end{proof}

\vskip 0.5cm Lemma \ref{principio-massimo} and (\ref{uguali}) imply that $\tilde W(\psi) = W(\psi)$, hence $\tilde E$ coincides with the energy $E$ for the functional $W$ on a neighbourhood of the point $(u,\omega)$. Moreover, since also $\tilde W'(\psi) = W'(\psi)$, then $(u,\omega)$ is a solution of (\ref{static}) with the potential $W$. Hence $(u,\omega)$ is a point of local minimum for the energy $E$ on the manifold $M_{C}$. Theorems \ref{thm-critical} and \ref{thm-stability} apply to $(u,\omega)$, and therefore $\psi(t,x)= u(x) e^{-i  \omega t}$ is a stable standing wave for (\ref{eq}).

This argument implies that we can assume condition (\ref{H3}) without loss of generality.

\vskip 0.5cm

\noindent \textbf{Proof of Theorem \ref{thm-critical}.} A point $
(u,\omega)\in Y$ is critical for $E(u,\omega)$ constrained on the manifold $
M_C$ if and only if there exists $\lambda \in $ such that
\begin{equation}  \label{Lagrvinc}
\begin{array}{c}
\frac{\partial E(u,\omega)}{\partial u}= \lambda\ \frac{\partial C(u,\omega)
}{\partial u} \\[0.3cm]
\frac{\partial E(u,\omega)}{\partial \omega}= \lambda\ \frac{\partial
C(u,\omega)}{\partial \omega}
\end{array}
\end{equation}
By definition of energy and charge (equations (\ref{energy-y}) and (\ref
{charge-y})), equation (\ref{Lagrvinc}) can be written as
\begin{equation*}
\begin{array}{rcl}
-\Delta u+ W^{\prime}(u)+\omega^2u & = & -2\lambda\, \omega\, u \\[0.3cm]
\int\omega\, u^2dx & = & -\lambda\int u^2dx.
\end{array}
\end{equation*}
Since $C\not= 0$, from the second equation $\lambda=-\omega$, and the first
becomes equation (\ref{static}). \qed

\vskip 0.5cm

\noindent \textbf{Proof of Lemma \ref{lemma-existence}.}  First of all we use (\ref{R-u}) to write the ratio $\Lambda(u,\omega)$ defined in (\ref{rapporto}) in the form
\begin{equation}  \label{lambda-J}
\Lambda(u,\omega)= \frac 1 2\, \left( \omega + \frac 1 \omega \right) +
\frac{J(u)}{|C(u,\omega)|}
\end{equation}
where
\begin{equation}  \label{funzionale-J}
J(u) := \int \left( \frac{1}{2} |\nabla u|^2 + R(u) \right)dx
\end{equation}

\noindent \emph{Step I.} If $(u_n,\omega_n)$ is a minimising sequence for $E$
on $M_{\bar C}$, then up to a subsequence
\begin{equation}  \label{converg}
\begin{array}{lll}
\omega_n \rightarrow \omega_0 \neq 0 &  &  \\
u_n \rightharpoonup u_0 \neq 0 \text{ weakly in } H^1(\R^N) &  &  \\
u_n \rightarrow u_0 \neq 0 \text{ strong in } L^t_{loc}(\R^N), \, 2\leq t<2^*.
&  &
\end{array}
\end{equation}
with
\begin{equation}  \label{rho}
\int u_n^2\, dx \to \rho:= \frac{|\bar C|}{\omega_0}
\end{equation}

\noindent To prove (\ref{converg}) we need to prove boundedness for $
(u_n,\omega_n)$. From (\ref{H0}) it follows that
\begin{equation}  \label{vicino-zero}
\exists\, \delta>0 \ \exists\, c_1>0, \text{ such that } W(z) \geq c_1
|z|^2 \text { for } 0\leq |z| \leq \delta.
\end{equation}
and one can choose $c_1=\min\limits_{0\leq |z| \leq \delta} \frac{W(z)}{\frac 1 2\, |z|^2}$
By equations (\ref{energy-y}) and (\ref{charge-y}) we can write
\begin{equation*}
E(u_n,\omega_n) =  \int \left( \frac{1}{2} |\nabla u_n|^2 + W(u_n)
\right) dx + \frac{|\bar C|\, \omega_n}{2}
\end{equation*}
Hence necessarily $\omega_n$ and $\int |\nabla u_n|^2 dx$ are bounded. Thus
up to a subsequence, $\omega_n \rightarrow \omega_0$. We now show that also $
\int u_n^2dx$ is bounded. Let us assume on the contrary that $(u_n,\omega_n)$
is minimising for $E$ on $M_{\bar C}$ and
\begin{equation*}
\int u_n^2\, dx \rightarrow \infty.
\end{equation*}
Clearly $\int W(u_n)\,dx$ is bounded and by (\ref{H0}) and (\ref{vicino-zero})
\begin{equation}  \label{for}
\int W(u_n)\, dx \geq \int_{0 \leq u_n\leq \delta} W(u_n)\, dx \geq c_1\,
\int_{0 \leq u_n \leq \delta} u_n^2\, dx .
\end{equation}
On the other hand
\begin{equation*}
\int_{0\leq u_n\leq \delta} u_n^2\, dx +\int_{u_n\geq \delta} u_n^2\, dx
\rightarrow \infty,
\end{equation*}
thus we have, by equation (\ref{for})
\begin{equation*}
\int_{u_n\geq \delta} u_n^2\, dx \rightarrow \infty.
\end{equation*}
This drives to a contradiction because
\begin{equation*}
\frac{1}{\delta^{2^*-2}} \ \int_{u_n\geq \delta} u_n^{2^*}dx \geq
\int_{u_n\geq \delta} u_n^2\, dx
\end{equation*}
and by the Sobolev embedding theorem
\begin{equation*}
\int_{u_n\geq \delta} u_n^{2^*}dx \leq \int u_n^{2^*}dx \leq K \int |\nabla
u_n|^2dx< const.
\end{equation*}
Thus, $\int|u_n|^2dx$ is bounded and the classical embedding theorems imply (\ref{converg}).

\vskip 0.3cm

\noindent \emph{Step II.} We prove that
\begin{equation}  \label{inf-uguali}
\inf\limits_{M_{\bar C}} E = \inf\limits_{N_\rho} J + \frac 1 2 \, \rho +
\frac{\bar C^2}{2\rho}
\end{equation}
where $\rho$ is defined in (\ref{rho}) and
\begin{equation*}
N_\rho:= \left\{u\in H^1(\R^N,\R^+)\ :\ \int u^2\, dx = \rho\right\}
\end{equation*}
By (\ref{funzionale-J}), the energy $E$ of any minimising sequence can be
written in the following form
\begin{equation*}
E(u_n,\omega_n)=J(u_n)+\frac{1}{2} \, \int u_n^2\, dx+\frac{\bar C^2}{2\int
u_n^2\, dx} \rightarrow m_c:=\inf\limits_{M_{\bar C}} E(u,\omega)
\end{equation*}
hence by (\ref{rho})
\begin{equation}  \label{inf-J}
J(u_n) \rightarrow m_c - \frac{1}{2}\, \rho - \frac{\bar C^2}{2\rho}.
\end{equation}
Now we take $v_n=\frac{\sqrt{\rho}}{||u_n||_{L^2}}u_n \in N_\rho$ and show
that
\begin{equation}  \label{convpalla}
J(v_n)-J(u_n)\rightarrow 0.
\end{equation}
Indeed,
{\small
\begin{equation*}
\begin{array}{l}
|J(v_n)-J(u_n)| = |(\frac{\rho}{||u_n||^2_{L^2}}-1) \int| \nabla u_n|^2dx +
\int R(\frac{\sqrt{\rho}}{||u_n||_{L^2}}u_n)-R(u_n)dx| \leq \\
\leq |(\frac{\rho}{||u_n||^2_{L^2}}-1)|\int|\nabla u_n|^2dx+ |(\frac{\sqrt{
\rho}}{||u_n||_{L^2}}-1) \int R^{\prime}(\theta\frac{\sqrt{\rho}}{
||u_n||_{L^2}}u_n+(1-\theta)u_n)u_n dx|
\end{array}
\end{equation*}}
and by (\ref{R-u}) and (\ref{H3})
{\small
\begin{equation*}
\begin{array}{l}
|\int R^{\prime}(\theta\frac{\sqrt{\rho}}{||u_n||_{L^2}}u_n+ (1-\theta) u_n)
u_n dx|\leq \\
c_1\int|\left(\theta\frac{\sqrt{\rho}}{||u_n||_{L^2}}+(1-\theta)
\right)^{q-1}||u_n^q| dx + c_2 \int |\left(\theta\frac{\sqrt{\beta}}{
||u_n||_{L^2}}+(1-\theta) \right)^{p-1}| |u_n^{p}|dx .
\end{array}
\end{equation*}}
The Sobolev embedding theorem gives that
\begin{equation*}
\int u_n^q\, dx, \ \ \int u_n^p\, dx \text { are bounded }
\end{equation*}
and we get equation (\ref{convpalla}) since $\frac{\rho}{||u_n||^2_{L^2}}
\rightarrow 1$. Hence by (\ref{inf-J}) and (\ref{convpalla}) it follows that
\begin{equation*}
J(v_n) \rightarrow m_c - \frac{1}{2}\, \rho-\frac{\bar C^2}{2\rho}.
\end{equation*}
hence
\begin{equation}  \label{inf-le}
\inf_{N_\rho} J \le m_c - \frac{1}{2}\, \rho-\frac{\bar C^2}{2\rho}.
\end{equation}

The opposite inequality holds since for any $v \in N_\rho$ the couple $
\left(v,\frac{|\bar C|}{\rho}\right) \in M_{\bar C}$ and
\begin{equation*}
E\left(v,\frac{|\bar C|}{\rho}\right)= J(v)+\frac{1}{2}\, \rho + \frac{\bar
C^2}{2\rho} \ge m_c
\end{equation*}
Hence
\begin{equation}  \label{inf-ge}
\inf_{N_{\rho}} J \ge m_c - \frac{1}{2}\, \rho + \frac{\bar C^2}{2\rho} .
\end{equation}
Equations (\ref{inf-le}) and (\ref{inf-ge}) together imply (\ref{inf-uguali}).

\vskip 0.3cm

\noindent \emph{Step III.} We now claim that the energy $E$ has a point of
minimum on the manifold $M_{\bar C}$ with $\bar C= C(\bar u, \bar \omega)$.

By hypothesis, there exists $(\bar u, \bar \omega) \in M_{\bar C}$ such that
$\Lambda(\bar u, \bar \omega) <1$. Hence by definition of $\Lambda$ (see (\ref{rapporto}))
\begin{equation}  \label{minimo-e-min-c}
\inf\limits_{M_{\bar C}}\, E < |\bar C|
\end{equation}
and moreover by elementary arguments
\begin{equation*}
\frac{1}{2}\, \rho + \frac{\bar C^2}{2\rho} \ge |\bar C|
\end{equation*}
for any $\rho>0$. Hence by (\ref{inf-uguali}) and (\ref{minimo-e-min-c}) it follows
\begin{equation}  \label{J-negativo}
\inf\limits_{N_\rho} J < 0
\end{equation}

Now by (\ref{J-negativo}) and by applying results in \cite{BBGM}, it follows
that the functional $J$ has a point of minimum $v_0$ constrained to the
manifold $N_\rho$, that is $J(v_0)=\inf_{N_\rho} J$. By repeating the same
argument as above and by (\ref{inf-uguali}) we have that
\begin{equation*}
E\left(v_0,\frac{|\bar C|}{\rho}\right)= J(v_0)+\frac{1}{2}\, \rho + \frac{
\bar C^2}{2\rho} = \inf\limits_{M_{\bar C}}\, E
\end{equation*}
and the thesis follows. \qed

\vskip 0.5cm The proof of Theorem \ref{main-theorem} is now completed by
proving Theorem \ref{thm-stability}.

\begin{rem}
\label{remark-conv} By applying results in \cite{BBGM}, we can say more
about any minimising sequence $(u_{n},\omega_{n})$ for $E$ on $M_{C}$. By
finding a correspondent minimising sequence for the functional $J$ of
equation (\ref{funzionale-J}) as done in step II, we obtain that the
sequence $(u_{n},\omega_{n})$ has the following properties: there exists a
sequence $(y_{n}) \subset \R^{N}$ such that $u_{n}(x)= \bar u(x+y_{n}) +
w_{n}(x)$ with $w_{n} \to 0$ in $H^{1}$; the limit $\bar \omega$ of the
sequence $\omega_{n}$ satisfies $\bar \omega^{2} < W^{\prime\prime}(0)$.
\end{rem}

\subsection{Proof of Theorem \ref{thm-stability}.} \label{proof-stability}

Let $C\in \R^{-}$ be such that the functional $E$ on $Y$ has a constrained
minimum on the manifold $M_{C}$, and all points of minimum $(u_{0},\omega_{0})$ in $M_{C}$ are isolated.

We first study relations between the energies $\E$ on $X$ and $E$ on $Y$
defined in Section \ref{wave-ds}. In (\ref{imm-y-x}) we have given the
embedding of $Y$ into $X$. If $(u,\omega) \in Y$ then $E(u,
\omega)=\E(u(x)e^{-i\omega t})$ and $C(u,\omega)=\Ch(u(x)e^{-i\omega t})$.
Let now $\mathbf{\Psi}=(\psi(t,x),\psi_{t}(t,x))$ be a function in $X$.
Using polar notation we write
\begin{equation*}
\psi(t,x)=u(t,x) e^{i S(t,x)}, \ u \ge 0
\end{equation*}
and set
\begin{equation}  \label{omega-tilde-psi}
\tilde \omega(t):= \frac{- \Ch(\psi)}{\int u^2(t,x) dx} = \frac{- \int
\partial_t S(t,x)\, u^2(t,x) dx}{\int u^2(t,x) dx}
\end{equation}
By using this notation we can associate to any function $\mathbf{\Psi} \in X$ a
point $(u(t,x),\tilde \omega(t))$ in $Y$ for any $t\in \R$. Hence we can
compare the energy $\E(\psi)$ of equation (\ref{energy}) and the energy $E(u,\tilde \omega)$
of equation (\ref{energy-y}).

\begin{lemma}
\label{due-energie} For any $\mathbf{\Psi}(t,x)$ in $X$ and for any $t \in \R$
it holds
\begin{equation}  \label{stima-ener}
\E(\mathbf{\Psi}) \ge E(u,\tilde \omega) + \frac 1 2 \int \, \left( |\nabla
S|^2 u^2 +|\partial_{t}u|^2 \right) dx
\end{equation}
\begin{equation}  \label{carica-n}
\Ch(\mathbf{\Psi}) = C(u,\tilde \omega)= - \tilde \omega\ \int u^2 dx
\end{equation}
\end{lemma}

\begin{proof}
Equation (\ref{stima-ener}) follows by applying H\"older
inequality to get
$$
\left(\int \partial_t S\, u^2 dx\right)^2 \le \int \partial_t
S^2\, u^2 dx \ \int u^2 dx
$$
and comparing (\ref{energy}) and (\ref{energy-y}). Equation (\ref{carica-n}) follows
from (\ref{omega-tilde-psi}).
\end{proof}

\begin{lemma}
\label{minimi-ener} For any $t\in \R$, it holds
\begin{equation*}
m_c=\inf\limits_{M_C}\ E(u,\omega)= \inf\limits_{\Ch(\mathbf{\Psi})=C}\
\E(\mathbf{\Psi})
\end{equation*}
\end{lemma}

\begin{proof}
>From Lemma \ref{due-energie}  it follows that if $\Ch(\mathbf{\Psi})=C$
then $C(u,\tilde \omega)=C$ and $\E(\mathbf{\Psi}) \ge E(u,\tilde \omega)$. Hence $\inf \E \ge \inf E(u,\tilde \omega) \ge \inf_{M_{C}} E$.

For the opposite inequality, we use the embedding (\ref{imm-y-x}). It implies that $\inf_{M_{C}} E = \inf_{M_{C}} \E(u(x)e^{-i\omega t}) \ge \inf \E$.
\end{proof}

It will be useful in the following to consider a kind of projection of a
couple $(u,\tilde \omega)\in Y$ on the manifold $M_{C}$. This is achieved by
defining for a function $\psi(t,x)$
\begin{equation}  \label{omega-psi}
\omega(t) := \frac{\tilde \omega(t) \, C}{\Ch(\mathbf{\Psi})}
\end{equation}
where $\tilde \omega(t)$ is as in (\ref{omega-tilde-psi}). Then it follows
that $C(u(t,x),\omega(t))= C$ for any $t\in \R$.

We now introduce the function
\begin{equation*}
V(\mathbf{\Psi})=\left( \E(\mathbf{\Psi})- m_c \right)^2+\left(\Ch(\mathbf{\Psi}) - C\right)^2
\end{equation*}
and prove that it is a Lyapunov function for the ground state set $\Gamma_C$
of equation (\ref{gamma-totale}) for any $t$. The first two conditions of
Definition \ref{lyapunov} follow straightforward from the definitions of
energy and charge, and from Lemma \ref{minimi-ener}. Indeed, if $V(\mathbf{
\Psi})=0$ then $\Ch(\mathbf{\Psi})=C$ and $\E(\mathbf{\Psi})=m_{c}$. Hence,
by Lemma \ref{due-energie} and \ref{minimi-ener}, it follows that $C(u,\tilde \omega)=C$ and $E(u,\tilde \omega)=m_{c}=\E(\mathbf{\Psi})$.
Hence $\nabla S = \partial_{t} u =0$ and $\partial_{t} S$ is constant almost
everywhere. Therefore $\mathbf{\Psi}$ is a standing wave and $(u,\tilde
\omega)$ is in $\Gamma_{C}$. Moreover, since energy and charge are first
integrals of the evolution flow of equation (\ref{eq-hamil}), $\frac{d}{dt}
V =0$. We now prove the last condition. We need to show that if $(\mathbf{\Psi_{n}})$ is a sequence of functions in $X$ satisfying $V(\mathbf{\Psi_{n}}
) \to 0$, then $d(\mathbf{\Psi_{n}}, \Gamma_{C}) \to 0$, where we recall
that the distance is given by the product topology on $X=H^{1}(\R^{N})\times
L^{2}(\R^{N})$.

Let $(\mathbf{\Psi_n}(t,x))$ be a sequence of functions in $X$ such that $\E(
\mathbf{\Psi_n}) \to m_c$ and $\Ch(\mathbf{\Psi_n})\to C$. Then using
notation introduced in (\ref{omega-tilde-psi}) and by Lemma \ref{due-energie}
and \ref{minimi-ener}, for any $t\in \R$ it holds $E(u_n,\tilde \omega_n) \to
m_c$, $C(u_n,\tilde \omega_n) \to C$ and
\begin{equation}  \label{stime-su-e-tu}
\int \left( |\nabla S_n|^2 u_n^2 +|\partial_{t}u_{n}|^2 \right) dx \to 0
\end{equation}
Hence, if $\omega_{n}(t)$ is defined as in (\ref{omega-psi}), we have
\begin{equation*}
E(u_{n}, \omega_{n}) \to m_{c}=\inf\limits_{M_{C}} E\, , \quad C(u_n,\tilde
\omega_n)=C \quad \forall\, t\in \R
\end{equation*}
that is $(u_{n},\omega_{n})$ is for any $t\in \R$ a minimising sequence for $E$
on $M_{C}$. By Remark \ref{remark-conv}, this implies that there exists a
sequence $(y_{n})$ in $\R^{N}$ such that up to a choice of sub-sequences
\begin{equation}  \label{prop-sottosuc}
u_{n}(x-y_{n}) \to u_{0}(x)\ \ \mbox{ in }\ \ H^{1} \qquad \mbox{and} \qquad
\omega_{n} \to \omega_{0}
\end{equation}
with $E(u_{0}, \omega_{0})=m_{c}$, $\omega_{0}<1$ and $(u_{0}, \omega_{0})$
is an isolated point of minimum for $E$ on $M_{C}$. Hence $\psi_{0}(t,x) =
u_{0}(x) e^{-i \omega_{0} t}$ is a standing wave solution of (\ref{eq}), and
to finish the proof we need to show that for any $t\in \R$ there exists $
\theta(t) \in \R$ such that
{\small
\begin{equation}  \label{aim-stab}
\| \psi_{n}(t,x-y_{n}) e^{-i \theta(t)} - \psi_{0}(t,x) \|_{H^{1}} + \|
\partial_{t} \psi_{n}(t,x-y_{n}) e^{-i \theta(t)} - \partial_{t}
\psi_{0}(t,x) \|_{L^{2}} \to 0
\end{equation}}
Notice that (\ref{aim-stab}) shows that the set $\Gamma(u_{0}, \omega_{0})$
(see (\ref{gamma-orbita})) is stable, hence the standing wave $\psi_{0}(t,x)
= u_{0}(x) e^{-i \omega_{0} t}$ is orbitally stable according to Definition
\ref{stabilita-orbitale}. However since $(u_{0}, \omega_{0})$ is an isolated
point of minimum then $\Gamma(u_{0}, \omega_{0})$ is a connected component
of $\Gamma_{C}$, hence the stability of $\Gamma_{C}$ is equivalent to that
of $\Gamma(u_{0}, \omega_{0})$.

\begin{lemma}
\label{norma-h1} For any $t\in \R$ there exists $\theta(t) \in \R$ such that
\begin{equation}  \label{aim-norma-h1}
\| \psi_{n}(t,x-y_{n}) e^{-i \theta(t)} - u_{0}(x) \|_{H^{1}} \to 0
\end{equation}
\end{lemma}

\begin{proof}
We first prove the lemma with a function $\theta(t,x) \in \R$, and later we show that $\theta(t,x)$ is almost everywhere constant in $x$. For clarity the dependence in $t$ is suppressed in formulas in which no confusion is possible.

By (\ref{prop-sottosuc}) and (\ref{stime-su-e-tu}) it follows that $u_{n}(x-y_{n})$ is bounded in $H^{1}$, and $\nabla(u_{n}(x-y_{n}))e^{iS_{n}(x-y_{n})}$ and $u_{n}(x-y_{n}) e^{iS_{n}(x-y_{n})} \nabla(S_{n}(x-y_{n}))$ are bounded in $L^{2}$, hence $\psi_{n}(x-y_{n})$ is bounded in $H^{1}$. By classical Sobolev embedding theorems it follows that up to a sub-sequence $\psi_{n}(x-y_{n})$ is convergent in $L^{2}_{loc}$ and almost everywhere to some function $\phi_{0}(t,x) \in L^{2}_{loc}$ for any $t\in \R$. We need to show that $\phi_{0}(t,x)= u_{0}(x) e^{i \theta(t,x)}$, where $u_{0}(x)$ is the radial component of the standing wave $\psi_{0}(t,x)$. By classical results on elliptic equations, we have that the function $u_{0}(x)$, which is a solution of (\ref{static}) with $\omega_{0}$, is smooth and positive. Hence
$$
\left| \frac{\psi_{n}(x-y_{n})}{u_{0}(x)} \right| \to \left| \frac{\phi_{0}(x)}{u_{0}(x)} \right| \quad \mbox{ almost everywhere}
$$
Moreover by (\ref{prop-sottosuc}), we have $u_{n}(x-y_{n}) \to  u_{0}(x)$ almost everywhere, hence
$$
\left| \frac{\phi_{0}(x)}{u_{0}(x)} \right| = 1 \quad \mbox{ almost everywhere}
$$
that is, for any $t\in \R$ there exists $\theta(t,x)$ such that
\begin{equation*}
\phi_{0}(t,x) = u_{0}(x) e^{i \theta(t,x)}
\end{equation*}
We recall that $\psi_{n}(x-y_{n})$ is bounded in $H^{1}$ and we have just shown that $\psi_{n}(x-y_{n}) \to u_{0}(x) e^{i \theta(t,x)}$ almost everywhere. Hence, using (\ref{prop-sottosuc}), it holds
{\small
$$
\psi_{n}(x-y_{n}) \to u_{0}(x) e^{i \theta(t,x)}\ \ \mbox{weakly in }\ \ L^{2}, \qquad \| \psi_{n}(x-y_{n}) \|_{L^{2}} \to \| u_{0}(x) e^{i \theta(t,x)} \|_{L^{2}}
$$}
hence
\begin{equation} \label{first-half}
\psi_{n}(x-y_{n}) \to u_{0}(x) e^{i \theta(t,x)}\ \ \mbox{ in }\ \ L^{2}
\end{equation}
To obtain (\ref{aim-norma-h1}) it now remains to prove that
\begin{equation} \label{second-half}
\nabla \psi_{n}(x-y_{n}) \to (\nabla u_{0}(x)) e^{i \theta(t,x)}\ \ \mbox{ in }\ \ L^{2}
\end{equation}
By (\ref{stime-su-e-tu}) $u_{n}(x-y_{n}) e^{iS_{n}(x-y_{n})} \nabla(S_{n}(x-y_{n}))$ vanishes for $n\to \infty$ in $L^{2}$, and $\nabla(u_{n}(x-y_{n}))e^{iS_{n}(x-y_{n})}$ is bounded in $L^{2}$. Moreover by (\ref{first-half}) it follows that $S_{n}(t,x) \to \theta(t,x)$ almost everywhere in $x$ for any $t$. Hence, using (\ref{prop-sottosuc}), it follows that
$$
\begin{array}{c}
\nabla(u_{n}(x-y_{n}))e^{iS_{n}(x-y_{n})} \to (\nabla u_{0}(x)) e^{i \theta(t,x)}\ \ \mbox{weakly in }\ \ L^{2}\\[0.3cm]
\qquad \| \nabla(u_{n}(x-y_{n}))e^{iS_{n}(x-y_{n})} \|_{L^{2}} \to \| (\nabla u_{0}(x)) e^{i\theta(t,x)} \|_{L^{2}}
\end{array}
$$
hence (\ref{second-half}) follows.

Now we prove that $\theta(t,x)$ is almost everywhere constant in $x$ for any $t$. We recall that $S_{n}(t,x) \to \theta(t,x)$ almost everywhere in $x$, and we can assume that $|S_{n}(t,x)| \le 2\pi$ for all $(t,x)$. Hence $S_{n} \to \theta$ in $L^{1}_{loc}$ and in the sense of distributions. Moreover, since $u_{n}(x-y_{n})/u_{0}(x) \to 1$ in $L^{2}_{loc}$, we conclude that
{\small
\begin{equation} \label{prima-dist}
\frac{u_{n}(x-y_{n})}{u_{0}(x)} \ \nabla S_{n}(t,x) \to \nabla \theta(t,x) \ \ \mbox{ in the sense of distributions}\ \ \forall \, t
\end{equation}}
Now, $u_{0}(x) >0$ for all $x$, hence for any compact set $B\subset \R^{N}$
{\small
$$
\int_{B}\, \left|  \frac{u_{n}(x-y_{n})}{u_{0}(x)} \ \nabla S_{n}(t,x) \right|\, dx \le \| u_{n}(x-y_{n}) \nabla(S_{n}(x-y_{n}))\|_{L^{2}} \, \left( \int_{B}\, \left|  \frac{1}{u_{0}(x)} \right|\, dx \right)^{\frac 1 2}
$$}
By (\ref{stime-su-e-tu}), it follows that
{\small
\begin{equation} \label{seconda-dist}
\frac{u_{n}(x-y_{n})}{u_{0}(x)} \ \nabla S_{n}(t,x) \to 0  \ \ \mbox{ in $L^{2}_{loc}$ and in the sense of distributions}\ \ \forall \, t
\end{equation}}
By (\ref{prima-dist}) and (\ref{seconda-dist}), it follows that $\nabla \theta (t,x)$ is a function in $L^{2}$ and  $\nabla \theta (t,x)=0$ almost everywhere.
\end{proof}

\begin{lemma}
\label{norma-l2} For any $t \in \R$ let $\theta(t) \in \R$ be chosen as in Lemma
\ref{norma-h1} . Then
\begin{equation*}
\| (\partial_{t} \psi_{n}(t,x-y_{n})) e^{-i \theta(t)} - (-i\omega_{0}
u_{0}(x)) \|_{L^{2}} \to 0
\end{equation*}
\end{lemma}

\begin{proof}
We write
{\small
$$
\partial_{t}   \psi_{n}(t,x-y_{n}) = \big((\partial_{t} u_{n}(t,x-y_{n})) + i  u_{n}(t,x-y_{n}) (\partial_{t} S_{n}(t,x-y_{n}))\big) e^{i S_{n}(t,x-y_{n})}
$$}
and the first term on the right vanishes for $n\to \infty$ in $L^{2}$ by (\ref{stime-su-e-tu}). Moreover, by the choice of the sequence $(\mathbf{\Psi_{n}}(t,x))$ and (\ref{omega-tilde-psi}), it holds $\E(\mathbf{\Psi_{n}}) - E(u_{n},\tilde \omega_{n}) \to 0$, which together with (\ref{stime-su-e-tu}) implies that for all $t\in \R$
{\small
\begin{equation} \label{uno-l2}
u_{n}(t,x-y_{n}) (\partial_{t} S_{n}(t,x-y_{n})) e^{i S_{n}(t,x-y_{n})} \to -\tilde \omega_{n}(t) u_{n}(t,x-y_{n}) e^{i S_{n}(t,x-y_{n})} \ \ \mbox{ in $L^{2}$}
\end{equation}}
Moreover by definition of $\omega_{n}(t)$ (see (\ref{omega-psi})) and by (\ref{prop-sottosuc}) it holds $\tilde \omega_{n} \to \omega_{0}$. Hence by Lemma \ref{norma-h1},
$$
\begin{array}{c}
\tilde \omega_{n}(t) u_{n}(t,x-y_{n}) e^{i S_{n}(t,x-y_{n})} \to \omega_{0} u_{0}(x) e^{i \theta(t)} \ \ \mbox{almost everywhere} \\[0.3cm]
\| \tilde \omega_{n}(t) u_{n}(t,x-y_{n}) e^{i S_{n}(t,x-y_{n})} \|_{L^{2}} \to \| \omega_{0} u_{0}(x) e^{i \theta(t)}  \|_{L^{2}}
\end{array}
$$
Repeating the same argument as in Lemma \ref{norma-h1} this implies that
\begin{equation} \label{due-l2}
\tilde \omega_{n}(t) u_{n}(t,x-y_{n}) e^{i S_{n}(t,x-y_{n})} \to \omega_{0} u_{0}(x) e^{i \theta(t)}  \ \ \mbox{ in $L^{2}$}\ \ \forall \, t
\end{equation}
By (\ref{uno-l2}) and (\ref{due-l2}) the lemma follows.
\end{proof}

Now (\ref{aim-stab}) follows by choosing $\theta(t) = -\omega_{0}t + \tilde
\theta(t)$, and the proof of Theorem \ref{thm-stability} is finished.

\subsection{Conclusive proofs} \label{proof-altre}

\noindent \textbf{Proof of Corollary \ref{123-carica}.} Let us consider the
family of functions $u_{R}\in H^{1}$ defined for $R\ge 0$ in (\ref{frittate}) and consider the frequencies $\omega_{R} = \sqrt{\alpha(u_{R})}$, where $
\alpha(u)$ is defined in (\ref{alpha}). Using (\ref{rapporto}) we find that
\begin{equation*}
\Lambda(u_{R},\omega_{R}) = \sqrt{\alpha(u_{R})}
\end{equation*}
and in Proposition \ref{alpha-0} we have shown that $\alpha(u_{R}) \to
\alpha_{0}$ as $R\to \infty$. By (\ref{H1}) we have $\alpha_{0}<1$, hence
there exists $R_{0}$ such that for all $R\ge R_{0}$ it holds $
\Lambda(u_{R},\omega_{R}) <1$. Hence let $C_{0} =
C(u_{R_{0}},\omega_{R_{0}}) $, then since
\begin{equation*}
C(u_{R},\omega_{R})= - \sqrt{\alpha(u_{R})}\, O(R^{N}) \to -\infty \quad
\mbox{as $R\to \infty$}
\end{equation*}
for any $C\in (-\infty, C_{0})$ there exists $R> R_{0}$ such that $
C(u_{R},\omega_{R})=C$ and $\Lambda(u_{R},\omega_{R}) <1$. Hence Lemma \ref
{lemma-existence} applies for existence, stability is obtained by Theorem \ref
{thm-stability} and the corollary is proved. \qed

\vskip 0.5cm

\noindent \textbf{Proof of Corollary \ref{123-carica-tutte}.} By using (\ref
{R-u}), (\ref{lambda-J}) and (\ref{funzionale-J}), we fix $\omega=1$ and we
have to prove that for any $C\in (-\infty,0)$ there exists $u\in H^{1}$ such
that $\Lambda(u,1)<1$ and $C(u,1)=C$. We let again $R\ge 0$ and define a
family of functions $u_{R}\in H^{1}$ by
\begin{equation}  \label{frittatissime}
u_{R}(x)= \left\{
\begin{array}{cl}
\gamma \, R^{-\frac N 2} & \mbox{ if $|x| \le R$} \\[0.3cm]
0 & \mbox{ if $|x| \ge 2R$} \\[0.3cm]
\gamma\, R^{-\frac N 2} \left( 2-\frac{|x|}{R} \right) &
\mbox{ if $R\le |x|
\le 2R$}
\end{array}
\right.
\end{equation}
where $\gamma$ is a constant to be fixed later. Moreover by (\ref{H1-prime})
we have for $R$ big enough
{\small
\begin{equation}  \label{stima-J-fritt}
J(u_{R}) = \int \left( \frac{1}{2} |\nabla u_{R}|^2 + W(u_{R}) - \frac 1 2
|u_{R}|^{2} \right)dx \le \int \left( \frac{1}{2} |\nabla u_{R}|^2 -
|u_{R}|^{2+\varepsilon} \right)dx
\end{equation}}
with $0<\varepsilon < \frac 4 N$. Hence (\ref{frittatissime}) and (\ref
{stima-J-fritt}) imply that
\begin{equation*}
J(u_{R}) \le const \, R^{-2}\, \left( 1 - R^{2- \varepsilon \frac N 2}
\right) \to 0^{-} \quad \mbox{as $R\to \infty$}
\end{equation*}
Hence let $R_{0}$ be such that $\Lambda(u_{R},1) <1$ for $R\ge R_{0}$. Then
we choose $\gamma$ such that $C(u_{R_{0}},1)=C$. The corollary is proved by
applying Lemma \ref{lemma-existence} for existence and Theorem \ref
{thm-stability} for stability. \qed


\begin{thebibliography}{99}

\bibitem{hylo}  J. Bellazzini, V. Benci, C. Bonanno, E. Sinibaldi, \emph{
Hylomorphic solitons}, in preparation

\bibitem{BBGM}  J. Bellazzini, V. Benci, M. Ghimenti, A.M. Micheletti, \emph{\
On the existence of the fundamental eigenvalue of an elliptic problem in $\R^N$
}, Adv. Nonlinear Stud. \textbf{7} (2007), 439--458

\bibitem{sammomme}  V. Benci, D. Fortunato, \emph{Solitary waves in the
nonlinear wave equation and in gauge theories}, J. Fixed Point Theory Appl.
\textbf{1} (2007), 61--86

\bibitem{Beres-Lions}  H. Berestycki, P.L. Lions, \textit{Nonlinear scalar
field equations. I. Existence of a ground state}, Arch. Rational Mech. Anal.
\textbf{82} (1982), 313--345

\bibitem{coleman78}  S. Coleman, V. Glaser, A. Martin, \emph{Action minima
among solutions to a class of euclidean scalar field equation}, Commun. Math.
Phys. \textbf{58} (1978), 211--221

\bibitem{Gelfand}  I.M. Gelfand, S.V. Fomin, ``Calculus of Variations",
Prentice-Hall, Englewood Cliffs, N.J. 1963

\bibitem{gss87}  M. Grillakis, J. Shatah, W. Strauss, \emph{Stability theory of
solitary waves in the presence of symmetry, I}, J. Funct. Anal. \textbf{74}
(1987), 160--197

\bibitem{Rajaraman}  R. Rajaraman, \textit{Solitons and instantons}, North
Holland, Amsterdam, Oxford, New York, Tokio, 1988

\bibitem{rosen68}  G. Rosen, \emph{Particle-like solutions to nonlinear
complex scalar field theories with positive-definite energy densities}, J.
Math. Phys. \textbf{9} (1968), 996--998

\bibitem{shatah}  J. Shatah, \emph{Stable standing waves of non-linear
Klein-Gordon equations}, Commun. Math. Phys. \textbf{91} (1983), 313--327

\bibitem{ss85}  J. Shatah, W. Strauss, \emph{Instability of nonlinear bound
states}, Commun. Math. Phys. \textbf{100} (1985), 173--190

\bibitem{strauss}  W.A. Strauss, \emph{Nonlinear invariant wave equations},
Lecture notes in physics, vol. 23, Springer, 1978

\bibitem{Witham}  G.B. Witham, \textit{Linear and nonlinear waves}, John
Wiley and Sons, New York, 1974
\end{thebibliography}
\end{document}